\documentclass{article}

\usepackage{amsmath}
\usepackage{array}
\usepackage{oldgerm}
\usepackage{times}
\usepackage{minitoc}
\usepackage{latexsym}
\usepackage{amsfonts}
\usepackage{psfrag}
\usepackage{amssymb,ifthen,amscd,amssymb}

\newtheorem{theo}{Theorem}[section]
\newtheorem{prop}{Proposition}[section]
\newtheorem{coro}{Corollary}[section]
\newtheorem{lemm}{Lemma}[section]
\newtheorem{defi}{Definition}[section]

\newtheorem{exam}{Example}
\newtheorem{rema}{Remark}[section]

\def\shuffle{\mathop{_{^{\sqcup\!\sqcup}}}} 
\gdef\stuffle{\;%
  \setlength{\unitlength}{0.0125cm}%
  \begin{picture}(20,10)(220,580)
  \thinlines
  \put(220,592){\line( 0,-1){ 10}}
  \put(220,582){\line( 1, 0){ 20}}
  \put(240,582){\line( 0, 1){ 10}}
  \put(230,592){\line( 0,-1){ 10}}
  \put(225,587){\line( 1, 0){ 10}}
  \end{picture}\;
}
\gdef\smuffle{\;%
  \setlength{\unitlength}{0.0125cm}%
  \begin{picture}(20,10)(220,580)
  \thinlines
  \put(220,592){\line( 0,-1){ 10}}
  \put(220,582){\line( 1, 0){ 20}}
  \put(240,582){\line( 0, 1){ 10}}
  \put(230,587){\circle*{5}}
  \end{picture}\;
}

\gdef\stufflemin{\;%
  \setlength{\unitlength}{0.0125cm}%
  \begin{picture}(20,10)(220,580)
  \thinlines
  \put(220,592){\line( 0,-1){ 10}}
  \put(220,582){\line( 1, 0){ 20}}
  \put(240,582){\line( 0, 1){ 10}}
  \put(225,587){\line( 1, 0){ 10}}
  \end{picture}\;
}
\gdef\duffle{\;%
  \setlength{\unitlength}{0.0175cm}%
  \begin{picture}(20,10)(220,580)
  \thinlines
  \put(220,592){\line( 0,-1){ 10}}
  \put(220,582){\line( 1, 0){ 20}}
  \put(240,582){\line( 0, 1){ 10}}
  \put(230,592){\line( 0,-1){ 10}}
  \put(222,587){\line( 1, 0){ 6}}
  \put(225,584){\line( 0, 1){ 6}}
  \put(235,587){\circle*{4}}   
  \end{picture}\;
}

\newcommand{\path}{{\rightsquigarrow}}

\newcommand{\C}{{\mathbb C}}
\newcommand{\N}{{\mathbb N}}

\newcommand{\Q}{{\mathbb Q}}
\newcommand{\Z}{{\mathbb Z}}

\newcommand{\bi}{\mathbf i}

\newcommand{\br}{\mathbf r}
\newcommand{\bs}{\mathbf s}
\newcommand{\bt}{\mathbf t}
\newcommand{\bF}{\mathbf F}

\newcommand{\bxi}{{\mathbf \xi}}
\newcommand{\brho}{{\mathbf \rho}}

\def\llm#1{\mathcal #1}

\newcommand{\calB}{\mathcal{B}}
\newcommand{\calC}{\mathcal{C}}
\newcommand{\calE}{\mathcal{E}}
\newcommand{\calI}{\mathcal{I}}

\newcommand{\calP}{\mathcal{P}}

\newcommand{\calT}{\mathcal{T}}

\newcommand{\ot}[1]{\overline{t}_{#1}}

\def\pol#1{\langle #1 \rangle}
\def\AX{A \pol{X}}

\newcommand{\serie}[2]{#1 \langle \! \langle #2 \rangle \! \rangle}
\def\AXX{\serie{A}{X}}

\newcommand{\Sum}[2]{\displaystyle{\sum_{#1}^{#2}}}

\newcommand{\ds}{\displaystyle}

\newcommand{\pg}{\geqslant}

\def\Di{\mathop{\rm Di}\nolimits}

\newcommand{\beginproof}{\noindent{\bf Proof}.~}
\def\endproof{%
  \unskip\nobreak\hfill\penalty80\hskip1em\hbox{}\nobreak\hfill$\Box$\\}

\begin{document}

\title{Combinatorial study of colored Hurwitz polyz\^etas}

\maketitle

\begin{center}
\large Jean-Yves Enjalbert$^{1}$, Hoang Ngoc Minh$^{1,2}$ 
\end{center}
\noindent{\em$1$ Université Paris 13, Sorbonne Paris Cité, LIPN, CNRS(, UMR 7030), F-93430, Villetaneuse, France.\\
$2$ Universit\'e Lille II,
1 place D\'eliot,
59024 Lille, France}\\[5pt]
 Email adresses : jean-yves.enjalbert@lipn.univ-paris13.fr, hoang@lipn.univ-paris13.fr  




\bigskip
\noindent{\bf abstract}\\
A combinatorial study discloses two surjective morphisms between generalized shuffle algebras and algebras generated by the colored Hurwitz  polyz\^etas.
The combinatorial aspects of the products and co-products involved in these algebras will be examined.



\section{Introduction}

\noindent
Classically, the Riemann z\^eta function is 
$\zeta(s)=\sum_{n>0}{n^{-s}},$
the Hurwitz z\^eta function is
$\zeta(s;t)=\sum_{n>0}{(n-t)^{-s}}$
and the colored z\^eta function is
$\zeta\binom{s}{q}=\sum_{n>0}{q^s}{n^{-s}},$
where $q$ is a root of unit. 
The three previous functions are defined over $\Z_{>0}$ but can be generalized over any {\em composition} (sequence of positive integers) $\bs=(s_1,\ldots,s_r)$, like, respectively, 
the  Riemann polyz\^eta function  
$\zeta(\bs)=\sum_{n_1>\ldots>n_r>0}{n_1^{-s_1}\ldots n_r^{-s_r}},$
the Hurwitz polyz\^eta function  
$\zeta(\bs;\bt)=\sum_{n_1>\ldots>n_r>0}{(n_1-t_1)^{-s_1}\ldots(n_r-t_r)^{-s_r}}$
and the colored polyz\^eta function
$\zeta\binom{\bs}{q^\bi}=\sum_{n_1>\ldots>n_r>0}{q^{i_1n_1}\ldots q^{i_rn_r}}{n_1^{-s_1}\ldots n_r^{-s_r}},$
with  $q$ a root of unit and $\bi=(i_1,\ldots,i_r)$ a composition.
These sums converge when $s_1>1$.

To study simultaneously these families of polyz\^etas, the colored Hurwitz  polyz\^etas,
for a composition $\bs=(s_1,\ldots,s_r)$ and a tuple of complex numbers $\bxi=(\xi_1,\ldots,\xi_r)$ and
a tuple of parameters in $]-\infty;1[$,  $\bt=(t_1,\ldots,t_r)$, are defined by \cite{HE10}
\begin{eqnarray}
\Di(\mathbf F_{\bxi,\bt};\bs)
=\sum_{n_1>\ldots>n_r>0}
\dfrac{\xi_1^{n_1}\ldots \xi_r^{n_r}}{(n_1-t_1)^{s_1}\ldots(n_r-t_r)^{s_r}}.
 \end{eqnarray}
Note that, for $l=1\ldots,r$, the numbers $\xi_l$ are not necessary roots of unity $q^{i_l}$.
We are working, in this note, with the condition\\
\centerline{$(E)\quad \forall i,\ \displaystyle|\prod_{k=1}^i\xi_k|\le 1\ \mbox{ and }\ t_i\in]-\infty;1[.$}
Hence, $\Di(\mathbf F_{\bxi,\bt};\bs)$ converges if $s_1>1$. We note $\calE$ the set of $\C$-tuples verifying $(E)$.

These  polyz\^etas are obtained as special values of iterated integrals\footnote{They are presented as generalized Nielsen polylogarithms in \cite{FPSAC95} (Definition 2.3) and as generalized Lerch functions in \cite{SLC44} (Definition 3).} over singular differential $1$-forms introduced in \cite{FPSAC95}. 
As  iterated integrals, they are encoded by words or by non commutative formal power series \cite{FPSAC95} and are used to construct bases for asymptotic expanding  \cite{cade} or symbolic integrating fuchian differential equations \cite{DM210} exactly or approximatively \cite{birkhauser}.
The meromorphic continuation of the colored Hurwitz polyz\^etas\footnote{See also references and a discussion about meromorphic continuation of Riemann polyz\^etas in \cite{HE07}.} is already studied in \cite{HE07,HE10}.
In our studies, we constructed an integral representation\footnote{This integral representation is obtained by applying successively the  polylogarithmic transform \cite{FPSAC95}. 
It is an application of {\em non} commutative convolution as shown in  \cite{hoang3} (Section 2.4). 
Other integral representations can be also deduced easily by change of variables,  for example $t=zr$ and then $r=e^{-u}$ \cite{HE07}.} of colored Hurwitz  polyz\^etas and  a {\em distribution} treating {\em simultanously two singularities} and our methods permit to  make the meromorphic continuation {\em commutatively} over the variables $s_1,\ldots,s_r$ \cite{HE07,HE10}. 
Moreover, \cite{HE10} gives another way to obtain the meromorphic continuation thanks to {\em translation equations} \cite{ecalle2003}. 
Our methods give the structure of multi-poles \cite{HE07} (Theorem 4.2) and two ways to calculate algorithmically the multi-residus\footnote{Other meromorphic continuations  can also be  obtained by Mellin transform as already done in \cite{Mat06} or by classical estimation on the imaginary part \cite{EMT11} but these later work reccursively, depth by depth, and  the commutativity of this process over the variables $s_1,\ldots,s_r$ must be proved. 
Unfortunately, the structure of multi-poles as well as multi-residus are missing in both works \cite{EMT11,Mat06}. 
In \cite{Kom10}, to make the meromorphic continuation (giving the expression of non positive integers multi-residus via a generalization of Bernoulli numbers -- but not of {\em all} multi-residus) of the specialization at roots of unity of colored Hurwitz polyz\^etas $\Di(\mathbf F_{\bxi,\bt};\bs)$,  the author bases on the  integral representation, on the contours, of the multiple Hurwitz-Lerch which corresponds {\em mutatis mutandis} to the  integral representation of generalized Lerch functions introduced earlier in \cite{HE07} (Corollary 3.3).}.

In this note, in continuation with our previous works \cite{FPSAC95,DM210,SLC44,words03,HE07,HE10},
we are focusing on {\em  Hofp algebra}, for a class of products as {\em minusstuffle} $(\stufflemin)$, {\em mulstuffle} $(\smuffle)$,  \ldots,
and in particular for the {\em  new} product {\em duffle} $(\duffle)$, obtained as ``tensorial product" of
$\smuffle$ and the well known stuffle $(\stuffle)$, of symbolic representations of these polyz\^etas
(see Definition \ref{Defstarproduct} and Proposition \ref{Prostarproduct} bellow).

\section{Combinatorial objects}

\subsection{Some products and their algebraic structures}

Let $X$ be an encoding alphabet and the free monoid over $X$ is denoted by $X^*$.
The {\em length} of any word $w\in X^*$ is denoted by  $|w|$ and
the unit of $X^*$ is denoted by $1_{X^*}$.
For any unitary commutative algebra $A$, a formal power series $S$ over $X$ with coefficients in $A$ can be written as the infinite sum
$\sum_{w\in X^*}\pol{S|w}w$.
The set of polynomials (resp. formal power series)  over $X$ with coefficients in $A$ is denoted by $\AX$ (resp. $\AXX$).
The set of degree $1$ monomials is  $AX=\{ax/ a\in A, x\in X\}$.

\begin{defi}\label{Defstarproduct}
We note $\calP$ the set of products $\star$ over $\AX$ verifying the conditions :
\begin{itemize}
\item[(i)] the map $\star : \AX\times\AX\rightarrow\AX$ is bilinear,
\item[(ii)] for any $w\in X^*$, $1_{X^*}\star w=w\star1_{X^*}=w$, 
\item[(iii)] for any $a, b\in X$ and $u,v\in X^*$, 
\begin{eqnarray*}
au\,\star\,bv 
= a(u\,\star\,bv)+b(au\,\star\,v)+[a,b](u\,\star\,v),
\end{eqnarray*}
where  $[.,.]:AX\times AX\rightarrow AX$ is a function verifying~: 
\begin{itemize}
\item[(S1)] $\forall a\in AX,\,[a,0]=0$ ,
\item[(S2)] $\forall (a,b)\in(AX)^2,\,[a,b]=[b,a]$,
\item[(S3)] $\forall (a,b,c)\in(AX)^3,\,[[a,b],c]=[a,[b,c]]$. 
\end{itemize}
\end{itemize}
\end{defi}
\begin{exam}[see \cite{reutenauer}]\label{ex:shuffle} {\em Product of  interated integrals}.\\
The shuffle is a bilinear product such that :
\begin{eqnarray*}
&\forall  w\in X^*&w\shuffle1_{X^*}=1_{X^*}\shuffle w=w
\qquad\mbox{ and}\cr
&\forall (a,b)\in X^2, \forall  (u,v)\in {X^*}^2, &au\shuffle vb 
= a(u\shuffle bv)+b(au\shuffle v)
.
\end{eqnarray*}
For example, for any letter $x_0$, $x$ and $x'$ in $X$,
\vskip-6mm 
\begin{eqnarray*}
x_0x'\shuffle x_0^2x = 
x_0x'x_0^2x+2x_0^2x'x_0x+3x_0^3x'x+3x_0^3xx'+x_0^2xx_0x'.
\end{eqnarray*}
\end{exam}
\begin{exam}[see \cite{Hof}]\label{ex:stuffle} {\em Product of quasi-symmetric functions}.\\
Let $X$ be an alphabet indexed by $\N$.\\
The stuffle is a bilinear product such that :
\begin{eqnarray*}
&&\forall  w\in X^*, 
\qquad w\stuffle1_{X^*}=1_{X^*}\stuffle w=w
\qquad
\mbox{ and}
\cr
&&\forall (x_i,x_j)\in X^2, \forall  (u,v)\in {X^*}^2,\cr
&&\qquad x_iu\stuffle x_jv 
= x_i(u\stuffle x_jv)+x_j(x_iu\stuffle v)+x_{i+j}(u\stuffle v)
.
\end{eqnarray*}
In particular, with the alphabet $Y=\{y_1,y_2,y_3,\ldots\}$,\vskip-7mm
\begin{eqnarray*}
(y_3y_1)\stuffle y_2
&=&y_3y_1y_2+y_3y_2y_1+y_3y_3+y_2y_3y_1+y_5y_1.
\end{eqnarray*}
\end{exam}
\begin{exam}[\cite{JSC}]\label{ex:minus-stuffle} {\em Product of  large multiple harmonic sums}.\\
Let $X$ be an alphabet indexed by $\N$.\\
The minus-stuffle is a bilinear product such that :
\begin{eqnarray*}
&&\forall  w\in X^*, 
\qquad w\stufflemin1_{X^*}=1_{X^*}\stufflemin w=w
\qquad
\mbox{ and}
\cr
&&\forall (x_i,x_j)\in X^2, \forall  (u,v)\in {X^*}^2,\cr
&&\qquad x_iu\stufflemin x_jv 
= x_i(u\stufflemin x_jv)+x_j(x_iu\stufflemin v)-x_{i+j}(u\stufflemin v)
.
\end{eqnarray*}
\end{exam}
\begin{exam}[\cite{HE10}]\label{ex:smuffle} {\em Product of colored sums}.\\
Let $X$ be an alphabet indexed by a monoid $(\calI,\times)$.\\
The mulstuffle is a bilinear product such that :
\begin{eqnarray*}
&&\forall  w\in X^*\qquad w\smuffle1_{X^*}=1_{X^*}\smuffle w=w
\qquad\mbox{ and}\cr
&&\forall (x_i,x_j)\in X^2, \forall  (u,v)\in {X^*}^2,\cr
&&\qquad x_iu\smuffle x_jv 
= x_i(u\smuffle x_jv)+x_j(x_iu\smuffle v)+x_{i\times j}(u\smuffle v)
.
\end{eqnarray*}
For example, with $X$ indexed by $\Q^*$,\\[2pt]
$x_{\frac23}x_{-1}\smuffle x_{\frac{1}{2}}
 = x_{\frac23}x_{-1}x_{\frac12}+x_{\frac23}x_{\frac12}x_{-1}
+x_{\frac23}x_{\frac{-1}{2}}
+x_{\frac{1}{2}}x_{\frac23}x_{-1}+x_{\frac13}x_{-1}.
$
\end{exam} 
\begin{rema}
Thanks to the one-to-one correspondence $(i_1,\ldots,i_r)\mapsto x_{i_1}\ldots x_{i_r}$ between tuples of $\calI$ and word over $X$, 
the calculus of $x_{\frac23}x_{-1}\smuffle x_{\frac{1}{2}}$ can be written as
$
\left(\frac23,-1\right)\smuffle\left(\frac{1}{2}\right)
=\left(\frac23,-1,\frac12\right)
+\left(\frac23,\frac12,-1\right)
+\left(\frac23,\frac{-1}{2}\right)
+\left(\frac{1}{2},\frac23,-1\right)
+\left(\frac13,-1\right).
$
\end{rema}
\begin{exam} [\cite{HE10}]{\em Product of colored Hurwitz polyz\^etas}.\\
Let $Y$ and $E$ be two alphabets and consider the alphabet $A=Y\times E$ with the concatenation defined recursively by $(y,e).(w_Y,w_E)=(yw_Y,ew_E)$ for any letters $y\in Y$, $e\in E$, and any word $w_Y\in Y^*$, $w_E\in E^*$.
The unit of the monoide $A^*$ is given by $1_{A^*}=(1_{Y^*},1_{E^*})$.
If $Y$ is indexed by $\N$ and $E$ by a monoid $(\calI,\times)$,
the duffle is a bilinear product such that
$\forall  w\in A^*,\quad w\duffle1_{A^*}=1_{A^*}\duffle w=w$,\\ 
$\forall (y_i,y_j)\in Y^2, \forall (e_l,e_k)\in E^2, \forall  (u,v)\in {A^*}^2,\quad
(y_i,e_l).u\duffle (y_j,e_k).v 
=(y_i,e_l).\allowbreak \left(u\duffle (y_j,e_k).v\right)
+(y_j,e_k).\left((y_i,e_l).u\duffle v\right)
+(y_{i+j},e_{l\times k}).(u\duffle v).$
\end{exam}
\begin{prop}\label{Prostarproduct}
The  shuffle, the  stuffle, the minus-stuffle and the mulstuffle  are elements of $\calP$, with respectively, $[x_i,x_j]=0,[x_i,x_j]=x_{i+j},[x_i,x_j]=-x_{i+j},[x_i,x_j]=x_{i\times j}$ for any letters $x_i$ and $x_j$ of $X$.\\
The duffle is in $\calP$, with $[(y_i,e_l),(y_j,e_k)]=(y_{i+j},e_{l\times k})$ for all 
$y_i, y_j$ in $Y$, $e_l, e_k$ in $E$.
\end{prop}
\begin{prop}
Let  $\star\in\calP$,
then $(\AX,\star)$ is a commutative algebra.
\end{prop}
\beginproof
We just have to show the commutativity and the associativity of $\star$.\\
To  obtain $w_1\star w_2=w_2\star w_1$ for all $w_1, w_2$ in $X^*$, we use an induction on $|w_1|+|w_2|$.
It is true when  $|w_1|+|w_2|\le 1$ thanks to (i) since $w_1$ or $w_2$ is $1_{X^*}$.
The equality (iii), the condition (S2) and the commutative of $+$ give the induction.
In the same way, an induction on  $|w_1|+|w_2|+|w_3|$ gives $w_1\star(w_2\star w_3)=(w_1\star w_2)\star w_3$
 thanks to (iii) and (S3).
\endproof
If we associate to each letter of $X$ an integer number called weight, the weight of a word is the sum of the weight of its letters.
In this case $X$ is graduated.\\
In \cite{Hof}, Hoffman works over  $\overline{X}=X\cup\{0\}$ with  $[.,.]:\overline{X}\times \overline{X}\rightarrow \overline{X}$ and call {\em quasi-product} any product in $\calP$ with the additional condition :
\begin{itemize}
\item[(S4)] Either $[a,b]=0$ for all $a, b$ in $X$; or the weight of $[a,b]$ is the sum of the weight of $a$ and the weight of $b$  for all $a, b$ in $X$.
\end{itemize}
\begin{exam}
\begin{enumerate}
\item The shuffle is a quasi-product.
\item Let $X$ be an alphabet indexed by $\N$ and define the weight of $x_i$, $i\in�\N$, by $i$ .
Then the stuffle is a quasi-product.
\end{enumerate}
\end{exam}
\begin{theo}[\cite{Hof}]\label{th:gradedHA}
If $X$ is graduated and has a quasi-product $\star$, then $(\AX,\star)$ is a commutative graduated $A$-algebra.. 
\end{theo}
\noindent\begin{tabular}[t]{lcl}
\hskip-6pt We can define&\qquad&(i) \ a comultiplication $\Delta : \AX\rightarrow\AX\otimes\AX$, \cr
&&(ii) \ a counit $\epsilon : \AX\rightarrow A$,
\end{tabular}\\
by : \quad $\forall w\in X^*$,\quad 
$\displaystyle 
\Delta w=\sum_{uv=w}u\otimes v$\quad and\quad$
\epsilon(w)=
\begin{cases}
1 \mbox { if }w=1_{X^*}\cr
0 \mbox{ otherwise}.
\end{cases}
$\\
The coproduct $\Delta$ is coassociative so $(\AX,\Delta,\epsilon)$ is a coalgebra.
\begin{lemm}\label{lem:structure}
For any $w\in X^*$ and $x\in X$,
$(x\otimes1_{X^*})\Delta w+1_{X^*}\otimes xw = \Delta xw.$
\end{lemm}
\beginproof
$\forall w\in X^*, \forall x\in X,
\displaystyle
\Delta xw
= \displaystyle\sum_{uv=xw}u\otimes v 
= \sum_{u'v=w}xu'\otimes v+1_{X^*}\otimes xw$\\
so\quad
$\displaystyle
\Delta xw
= \displaystyle x\otimes1_{X^*}\Big(\sum_{u'v=w}u'\otimes v\Big)+1_{X^*}\otimes xw 
= (x\otimes1_{X^*})\Delta w+1_{X^*}\otimes xw.
$
\endproof
\goodbreak
\begin{prop}
If $\star\in\calP$, then $(\AX,\star,\Delta,\epsilon)$ is a bialgebra.
\end{prop}
Remember that $\star$ acts over $\AX\otimes\AX$ by $(u\otimes v)\star(u'\otimes v')=(u\star u')\otimes(v\star v')$.
\goodbreak
\beginproof
$\epsilon$ is obviously a $\star$-homomorphism. 
It still has to be show 
$\Delta(w_1)\star\Delta(w_2)=\Delta(w_1\star w_2) $ over $X^*$.
This equality is true if $w_1$ or $w_2$ is equal to $1_{X^*}$.\\
Assume now that  $\Delta(u)\star\Delta(v)=\Delta(u\star v)$ for any word $u$ and $v$ such that $|u|+|v|\le n$, $n\in\N$, and let $w_1$ and $w_2$ be in $X^*$ with $|w_1|+|w_2|=n+1$.
We note $w_1=au$ and $w_2=bv$, with $a$ and $b$ two letters of $X$, $u$ and $v$ two words of $X^*$.
Thus, by definition,
$\Delta w_1=\sum_{u_1u_2=u}au_1\otimes u_2+1_{X^*}\otimes au$
and
$\Delta w_2=\sum_{v_1v_2=v}bv_1\otimes v_2+1_{X^*}\otimes bv$.
\begin{eqnarray*} 
&&\Delta(w_1)\star\Delta(w_2)\\
&&=
\sum_{{u_1u_2=u},{v_1v_2=v}}(au_1\star bv_1)\otimes(u_2\star v_2)
+\sum_{u_1u_2=u}(au_1)\otimes(u_2\star bv)\cr
&&+\sum_{v_1v_2=v}(bv_1)\otimes(au\star v_2)
+1_{X^*}\otimes(au\star bv)\cr
&&=
\sum_{{u_1u_2=u},{v_1v_2=v}}
(a(u_1\star bv_1)\otimes(u_2\star v_2)
+b(au_1\star v_1)\otimes(u_2\star v_2)\cr
&&+([a,b](u_1\star v_1))\otimes(u_2\star v_2))
+\sum_{u_1u_2=u}(au_1)\otimes(u_2\star bv)\cr
&&+\sum_{v_1v_2=v}(bv_1)\otimes(au\star v_2)
+1_{X^*}\otimes a(u\star bv)\\
&&+1_{X^*}\otimes b(au\star v)+1_{X^*}\otimes [a,b](u\star v)\\
&&=
\sum_{{u_1u_2=u},{v_1v_2=v}}a(u_1\star bv_1)\otimes(u_2\star v_2)+\sum_{u_1u_2=u}(au_1)\otimes(u_2\star bv)\cr
&&
+\sum_{{u_1u_2=u},{v_1v_2=v}}b(au_1\star v_1)\otimes(u_2\star v_2)+\sum_{v_1v_2=v}(bv_1)\otimes(au\star v_2)\cr
&&+
[a,b]\otimes1_{X^*}\sum_{\stackrel{u_1u_2=u}{v_1v_2=v}}(u_1\otimes u_2)\star (v_1\otimes v_2)\cr
&&
+(1_{X^*}\otimes a(u\star bv)+1_{X^*}\otimes b(au\star v)+1_{X^*}\otimes[a,b](u\star v))\\
&&=(a\otimes1_{X^*})(\Delta(u)\star\Delta(w_2))+1_{X^*}\otimes a(u\star bv)
+(b\otimes1_{X^*})(\Delta(w_1)\star\Delta(v))\cr
&&
+1_{X^*}\otimes b(au\star v)
+([a,b]\otimes1_{X^*})(\Delta(u)\star\Delta(v))+1_{X^*}\otimes [a,b](u\star v).
\end{eqnarray*} 
Using the induction hypothesis then the lemma \ref{lem:structure} (since $[a,b]\in AX$) gives
\begin{eqnarray*} 
\Delta(w_1)\star\Delta(w_2)
&=&\Delta(a(u\star w_2))+\Delta(b(w_1\star v))+\Delta([a,b](u\star v))\cr
&=&\Delta(a(u\star w_2)+b(w_1\star v)+[a,b](u\star v))\cr
&=&\Delta(w_1\star w_2).
\end{eqnarray*}
\endproof
\begin{rema}
In particular,
$\Delta$ is a~$\shuffle$-homomorphism, a~$\stuffle$-homomorphism and a~$\smuffle$-homomorphism.
\end{rema}
Let $\calC_n$ be the set of positive integer sequences $(i_1,\ldots,i_k)$ such that 
$i_1+\ldots+i_k\!=\!n$.
\begin{theo} Define $a_\star$ by, for all $x_1,\ldots,x_n$ in $X$,
\begin{eqnarray*}
&&a_\star(x_1\ldots x_n)\cr 
&&=\sum_{(i_1,\ldots,i_k)\in\calC_n}
(-1)^kx_1\ldots x_{i_1}\star x_{i_1+1}\ldots x_{i_1+i_2}\star\ldots\star x_{i_1+\ldots+i_{k-1}+1}\ldots x_{n}
\end{eqnarray*}
then, if $\star\in\calP$,  $(\AX,\star,\Delta,\epsilon,a_\star)$ is a Hopf algebra.
\end{theo}
\beginproof
With the applications :\\
\begin{tabular}{llcl}
$\mu $ :&$A$&$\rightarrow$&$\AX$\cr
&$\lambda$&$\mapsto$&$\lambda\,1_{X^*}$
\end{tabular}\quad
and \quad
\begin{tabular}{llcl}
$m $ :&$\AX\otimes\AX$&$\rightarrow$&$\AX$\cr
&$u\otimes v$&$\mapsto$&$u\star v$
\end{tabular},\\
the antipode must verify
$
m\circ(a_\star\otimes Id)\circ\Delta=\mu\circ\epsilon
$, or, in equivalent terms 
\begin{eqnarray*}
\sum_{uv=w}a_\star(u)\star v=\langle w|1_{X^*}\rangle1_{X^*}.
\end{eqnarray*}
i.e.
$\begin{cases}
a_\star(1_{X^*})=1_{X^*}n\cr 
\forall x\in X, a_\star(x)=-x\   
\end{cases}$ 
and, if $w=x_1\ldots x_n$ with $n\ge2$, $x_1,\ldots, x_n\in X$, 
\begin{eqnarray*}
a_\star(w)=-\sum_{k=1}^{n-1}a_\star(x_1\ldots x_k)\star x_{k+1}\ldots x_n.
\end{eqnarray*}
An induction over the length $n$ shows that $a_\star$ defined in theorem verifies these equalities, and, in the same way, $a_\star$ verifies 
$
m\circ(Id\otimes a_\star)\circ\Delta=\mu\circ\epsilon
$.
\endproof
\begin{coro}
If $\star$ is $\shuffle$ or $\stuffle$ or $\smuffle$ or $\duffle$, then this construction gives an Hopf algebra.
Moreover, for  $\shuffle$ or $\stuffle$, we obtain a graduated Hopf algebra.
\end{coro}

\subsection{Iterated integral}

Let us associate to each letter $x_i$ in $X$ a $1$-differential
form $\omega_i$, defined in some connected open subset $\llm U$ of
$\C$. 
For all paths $z_0\path z$ in $\llm U$, the {\em Chen
iterated integral} associated to $w=x_{i_1}\cdots x_{i_k}$ along
$z_0\path z$, noted  is defined recursively as follows
\begin{eqnarray}\label{intint}
\alpha_{z_0}^z(w)=\int\limits_{z_0\path z}\omega_{i_1}(z_1)\alpha_{z_0}^{z_1}(x_{i_2}\cdots x_{i_k})
&\mbox{and}&\alpha_{z_0}^z(1_{X^*})=1,
\end{eqnarray}
verifying the {\em rule of integration by parts} \cite{Chen} : 
\begin{eqnarray}
\alpha_{z_0}^z(u\shuffle
v)=\alpha_{z_0}^z(u)\alpha_{z_0}^z(v).
\end{eqnarray}
We extended this definition over $\AX$ (resp. $\AXX$) by 
\begin{eqnarray}
\displaystyle \alpha_{z_0}^z(S) = \sum_{w\in X^*}\pol{S|w}\alpha_{z_0}^{z}(w).
\end{eqnarray}

\goodbreak

\subsection{Shuffle relations}

\subsubsection{First encoding for colored Hurwitz polyz\^etas}

Let $\bxi=(\xi_n)$ be a sequence of complex numbers and $T$ a family of parameters.
Put $X'$ an alphabet indexed over $\N^*\times\C^\N\times T$ and $X=\{x_0\}\cup X'$.  
To each $x$ in $X$ we associate the differential form :
\begin{eqnarray}
\begin{cases}
\omega_{0}(z)=\dfrac{dz}{z} &\mbox{ si }\,x=x_0\cr
\omega_{i,\bxi,t}(z)=\dfrac{\prod_{k=1}^i\xi_k}{1-\prod_{k=1}^i\xi_k\,z}\times \dfrac{dz}{z^{t}}  &\mbox{ if }\,x=x_{i,\bxi,t}\mbox{ with }i\pg1.
\end{cases}
\end{eqnarray}
For any $T$-tuple $\bt=(t_1,\ldots,t_r)$ we associate the $T$-tuple $\overline{\bt}=(\overline{t_1},\ldots,\overline{t_r})$ given by
\begin{eqnarray}
\begin{cases}
\overline{t_1}&=\ t_1-t_2,\cr
\overline{t_2}&=\ t_2-t_3,\cr
&\vdots\cr
\overline{t_r}&=\ t_{r-1}-t_r
\end{cases}
\quad\mbox{in this way}\quad
\begin{cases}
t_1&=\ \overline{t_1}+\overline{t_2}+\ldots+\overline{t_n},\cr
t_2&=\ \overline{t_2}+\ldots+\overline{t_n},\cr
&\vdots\cr
t_r&=\ \overline{t_r}
\end{cases}
\end{eqnarray}
We choose the sequence $\xi$ and the family $\bt$ such that the condition $(E)$ is satisfied.
\begin{prop}\label{pro:1code}
For any $\bs=(s_1,\ldots,s_r)$ with $s_1>1$
if $\xi=(\xi_1,\ldots,\xi_r)\in\calE^r$ and $\bt=(t_1,\ldots,t_r)\in T^r$, then
$\Di(\mathbf F_{\bxi,\bt};\bs)=\alpha_0^1(x_0^{s_1-1}x_{1,\bxi,\overline{t_1}}\ldots x_{0}^{s_r-1}x_{r,\bxi,\overline{t_r}}).$
\end{prop}
\beginproof 
Since
$\displaystyle
\omega_{i,\bxi,t}(z) =\sum_{n>0}\prod_{k=1}^i\xi_k^n\dfrac{z^ndz}{z^{1+t}}$
then $ \displaystyle\alpha_0^z(x_{r,\bxi,\overline{t_r}})=\sum_{n>0}\prod_{k=1}^r\xi_k^n\dfrac{z^{n-\overline{t_r}}}{n-\overline{t_r}}$ and
$\displaystyle\alpha_0^z(x_0^{s_r-1}x_{r,\bxi,\overline{t_r}})=\sum_{n>0}\prod_{k=1}^r\xi_k^n\dfrac{z^{n-\overline{t_r}}}{(n-\overline{t_r})^{s_r}}$.
Hence, $\alpha_0^1(x_0^{s_1-1}x_{1,\bxi,\overline{t_1}}\ldots x_{0}^{s_r-1}x_{r,\bxi,\overline{t_r}})$ gives
$\displaystyle 
\sum_{m_1,\ldots,m_r>0}\prod_{j=1}^r\dfrac{\prod_{k_j=1}^j\xi_{k_j}^{m_j}}{{(m_j+\ldots+m_r-\ot{j}-\ldots-\ot{r})}^{s_j}}$, and then, by change of variables,
$\displaystyle\sum_{n_1>\ldots>n_r>0}\dfrac{\xi_1^{n_1}\ldots\xi_r^{n_r}}{(n_1-t_1)^{s_1}\ldots(n_r-t_r)^{s_r}}$.
\endproof
\vskip-5mm
\begin{theo}
Let $\calT$ be the group of parameters generated by $\langle T;+\rangle$, $\calC$ be a sub-group of $(\C^*,.)$ and $A$ a sub-ring of $\C$.
Put $\calC'=\calC^\N\cap\calE$ and $\calT'$ the set of finite tuple with elements in $\calT$.
Then the $A$ algebra generated by $\{\Di(\mathbf F_{\bxi,\bt};\bs)\}_{\bxi\in\calC', \bt\in\calT'}$ is the $A$ modulus  generated by $\{\Di(\mathbf F_{\bxi,\bt};\bs)\}_{\bxi\in\calC', \bt\in\calT'}$.
\end{theo}
\beginproof
We have express the product $\Di(\mathbf F_{\bxi,\bt};\bs)\Di(\mathbf F_{\bxi',\bt'};\bs')$, with $\bs=(s_1,\ldots,s_r)$,
$\bs'=(s_1',\ldots,s_{r'}')$, $\bxi,\bxi'\in\calC'$ and $\bt=(t_1,\ldots,t_r),\bt'=(t_1',\ldots,t_r')\in\calT'$, as linear combination of colored Hurwitz  polyz\^etas.
This is an iterated integral  associated to $x_0^{s_1-1}x_{1,\bxi,\overline{t_1}}\ldots x_{0}^{s_r-1}x_{r,\bxi,\overline{t_r}}\shuffle x_0^{s_1'-1}x_{1,\bxi',\overline{t_1'}}\ldots x_{0}^{s_{r'}'-1}x_{r',\bxi',\overline{t_{r'}'}}$
which is a sum of terms of the form 
$x_0^{s_1''-1}x_{1,\xi^{(1)},\overline{t_{(1)}}}\ldots x_0^{s_i''-1}x_{j_i,\xi^{(i)},\overline{t_{(i)}}}\ldots x_0^{s_r''-1}x_{j_{r''},\xi^{(r'')},\overline{t_{(r'')}}}$,
with $s_i''\in\N$, $\xi^{(i)}$ is $\bxi$ or $\bxi'$ and $t_{(i)}$ is $t_{j_i}$ or $t_{j_i}'$  for all $i$; and $r''=r+r''$.
Note that
\begin {eqnarray*}
&&\alpha_0^z(x_{i,\bxi,\overline{t_i}}x_0^{s-1}x_{j,\bxi',\overline{t_j}})\cr
&&=\int_0^z\sum_{m>0}\prod_{k=1}^i\xi_k^mz_1^{m-\overline{t_i}-1}dz_1\int_0^{z_1}\dfrac{dz_2}{z_2}...\int_0^{z_{s+1}}\sum_{n>0}\prod_{k=1}^i\xi_k'^nz_{s+1}^{n-\overline{t_j'}-1}dz_{s+1}\cr
&&=\sum_{m,n>0}\dfrac{\left(\xi_1\ldots\xi_i\right)^m\left(\xi_1'\ldots\xi_j'\right)^n}{(m+n-\overline{t_i}-\overline{t_j'})(n-\overline{t_j'})^{s}}\,z^{n+m},\\
&&\alpha_0^1\big(x_0^{s_1''-1}x_{1,\xi^{(1)},\overline{t_{(1)}}}\ldots x_0^{s_i''-1}x_{j_i,\xi^{(i)},\overline{t_{(i)}}}\ldots x_0^{s_r''-1}x_{j_{r''},\xi^{(r'')},\overline{t_{(r'')}}}\big)\qquad\cr
&&=\sum_{m_1,\ldots,m_{r''}>0}\prod_{i=1}^{r''}\dfrac{(\xi_1^{(i)}\ldots\xi_{j_i}^{(i)})^{m_i}}{(m_i+\ldots+m_{r''}-\overline{t_{(i)}}-\ldots-\overline{t_{(r'')}})^{s_i''}}\cr
&&=\sum_{n_1>\ldots>n_{r''}>0}\prod_{i=1}^{r''}\dfrac{\xi_i''^{\,n_i}}{(n_i-t_{i}'')^{s_i''}}
\end{eqnarray*}
with $n_i=m_i+\ldots+m_{r''}$, $t_i''=\overline{t_{(i)}}+\ldots+\overline{t_{(r'')}}$ for all $i$, so $\bt''\in\calT$;
$\xi_1''=\xi_1^{(1)}$ and $\xi_i''=\frac{\xi_1^{(i)}\ldots\xi_{j_i}^{(i)}}{\xi_1^{(i-1)}\ldots\xi_{j_{i-1}}^{(i-1)}}$ for $i>1$ so $\bxi''\in\calC$ : we can express each term of the shuffle product as $\Di(\mathbf F_{\bxi'',\bt''};\bs'')$.
\endproof
Note that the shuffle product over two words of  $X^*X'$ acts separately over $(\calC',.)$, $(\calT',+)$ and the convergent compositions.
We can describe the situation with the shuffle algebra\footnote{Working in $\Q\langle\left(x_0^*x_{i,\xi,t}\right)^*\rangle$ implies working in the graduated Hopf algebra $(\Q\langle X^*\rangle,\shuffle,\Delta,\epsilon,a_{\shuffle})$.}~:
\begin{theo}\label{th:Hcp}
Let $\cal H$ be the $\Q$-algebra generated by the colored Hurwitz  polyz\^etas. The map
$\zeta: (\Q\langle\left(x_0^*x_{i,\xi,t}\right)^*\rangle,\shuffle)\twoheadrightarrow({\cal H},.)$,
$x_0^{s_1}x_{1,\bxi,\overline{t_1}}\ldots x_{0}^{s_r}x_{r,\bxi,\overline{t_r}}\mapsto\Di(\mathbf F_{\bxi,\bt};\bs+1)$
 is a surjective algebra morphism.
\end{theo}
\begin{exam}
Since
$\Di(\mathbf F_{\xi,t};3)=\alpha_0^1(x_0^2x_{1,\xi,t})$  and $\Di(\mathbf F_{\xi',t'};2)=\alpha_0^1(x_0x_{1,\xi',t'})$
then $\Di(\mathbf F_{xi,t};3)\Di(\mathbf F_{xi',t'};2)= \alpha_0^1(x_0x_{1,\xi',t'}\shuffle x_0^2x_{1,\xi,t}).$
Example \ref{ex:shuffle} with $x=x_{1,\xi,t}$ and $x'=x_{1,\xi',t'}$ gives the expression of
$x_0x_{1,\xi',t'}\shuffle x_0^2x_{1,\xi,t}$.
But the first term obtained is
\begin{eqnarray*}
&&\alpha_0^1(x_0x_{1,\xi',t'}x_0^2x_{1,\xi,t})\cr
&&=
\int_0^1\dfrac{dz_1}{z_1}\int_0^{z_1}\sum_{m>0}{\xi'\,}^mz_2^{m-t'-1}dz_2
\int_0^{z_2}\dfrac{dz_3}{z_3}\int_0^{z_3}\dfrac{dz_4}{z_4}
\int_0^{z_4}\sum_{n>0}\xi^nz_5^{n-t-1}dz_5\cr
&&=\sum_{n,m>0}\dfrac{{\xi'\,}^m\xi^n}{(m+n-t'-t)^2(n-t)^3}\cr
&&=\sum_{n_1>n_2>0}\dfrac{{(\xi')}^{n_1}({\xi}/{\xi'})^{n_2}}{(n_1-t'-t)^2(n_2-t)^3}\\[2pt]
&&=\Di(\mathbf F_{(\xi,\xi/\xi');(t+t',t)};(2,3)).
\end{eqnarray*}

We can make similar calculus for the other terms and  find :
\begin{eqnarray*}
&&\Di(\mathbf F_{\xi,t};3)\Di(\mathbf F_{\xi',t'};2)\cr
&&=
\Di(\mathbf F_{(\xi',\xi/\xi');(t+t',t)};(2,3))
+2\Di(\mathbf F_{(\xi',\xi/\xi');(t+t',t)};(3,2))\cr
&&\ \ +\,3\Di(\mathbf F_{(\xi',\xi/\xi');(t+t',t)};(4,1))
+3\Di(\mathbf F_{(\xi,\xi'/\xi);(t+t',t')};(4,1))\cr
&&\ \ +\Di(\mathbf F_{(\xi,\xi'/\xi);(t+t',t')};(3,2)).
\end{eqnarray*}
\end{exam}

\subsubsection{Second encoding for colored Hurwitz polyz\^etas}

For the Hurwitz polyz\^etas, we can obtain an encoding indexed  by a finite alphabet.
Let the alphabet $X=\{x_0;x_1\}$ and associate to $x_0$ the form $\omega_0(z)=z^{-1}{dz}$ and at $x_1$ the form $\omega_1(z)=(1-z)^{-1}{dz}$.

For each $x\in X$ and $\lambda\in\C$, we note $(\lambda x)^*=\sum_{k\ge0}(\lambda x)^k$.
Then, (see \cite{FPSAC95}, \cite{DM210}),
$\alpha_0^1\left(x_0^{s_1-1}(t_1x_0)^{*s_1}x_1\ldots x_0^{s_r-1}(t_rx_0)^{*s_r}x_1\right) = \zeta(\bs;\bt).$

\begin{theo}\label{th:Hpc}
Let $\cal H'$ be the $\Q$-algebra generated by the Hurwitz polyz\^etas and $\cal X$ the $\Q$-algebra generated by $(t_1x_0)^{*s_1}x_1\ldots (t_rx_0)^{*s_r}x_r$.
Then, 
$\zeta:(\cal X,\shuffle)\twoheadrightarrow(\cal H',.)$ is a surjective morphism of algebras. 
\end{theo}

Note that we can apply the idea of encoding of ``simple'' colored Hurwitz  zetas functions (with depth one : $r=1$).
Let  $\bxi=(\xi_n)$ be a sequence of complex numbers in the unit ball $\calB(0;1)$ and $T$ a family of parameters.
Let $X=\{x_0,x_1,\ldots\}$ be a alphabet indexed by $\N$.
Associate to $x_0$ the differential form $\omega_0(z)=z^{-1}{dz}$ and to $x_i$, $i\ge1$, the differential form $\omega_i(z)=\xi_i(1-\xi_iz)^{-1}dz$.
\begin{prop}
With this notation,
$\alpha_0^1\left(\left((tx_0)^*x_0\right)^{s-1}(tx_0)^*x_i\right)=\Sum{n>0}{}\dfrac{\xi_i^n}{(n-t)^s}.$
\end{prop}
\beginproof Since $\ds\dfrac{\xi_idz_0}{1-\xi_iz_0}=\xi_i\sum_{n\ge0}(\xi_iz_0)^ndz_0$, we can write
\begin{eqnarray*}
\alpha_0^z\left((tx_0)^kx_i\right)
 = t^k\int_0^z\dfrac{dz_k}{z_k}\int_0^{z_k}\ldots\int_0^{z_1}\xi_i\sum_{n\ge0}(\xi_iz_0)^ndz_0 
= \sum_{n>0}t^k\dfrac{\xi_i^nz^n}{n^{k+1}},
\end{eqnarray*}
for $z\in \calB(0;1)$ and for $k\in\N$. 
Thanks to the absolute convergence,
\begin{eqnarray*}
\alpha_0^z\left((tx_0)^*x_i\right)
=\sum_{n>0}\dfrac{\xi_i^nz^n}{n}\sum_{k\ge0}\left(\dfrac{t}{n}\right)^k
=\sum_{n>0}\dfrac{\xi_i^nz^n}{n-t}.
\end{eqnarray*}
In the same way, if $z\in\calB(0;1)$ :
\begin{eqnarray*}
&\forall k\in\N,&\quad \alpha_0^z\left((tx_0)^kx_0(tx_0)^*x_i\right)=\Sum{n>0}{}t^k\dfrac{\xi_i^n}{n-t}\dfrac{z^n}{n^{k+1}},\cr 
&\mbox{so}\quad
&\alpha_0^z\left((tx_0)^*x_0(tx_0)^*x_i\right)=\Sum{n>0}{}\dfrac{\xi_i^nz^n}{(n-t)^2}\cr 
&\mbox{and}\quad&
\alpha_0^z\left(\left((tx_0)^*x_0\right)^{s-1}(tx_0)^*x_i\right)=\Sum{n>0}{}\dfrac{\xi_i^nz^n}{(n-t)^s}.
\end{eqnarray*}
\endproof
\begin{rema}
Note that, with the same notation,
\begin{eqnarray*}
\alpha_0^z\left(x_1\left((t_2x_0)^*x_0\right)^{s-1}(t_2x_0)^*x_2\right)
&=&\Sum{n,m>0}{}\dfrac{\xi_2^n\xi_1^mz^{n+m}}{(n-t_2)^s(m+n)}\cr
&=&\Sum{n_1>n_2>0}{}\dfrac{\xi_2^{n_2}\xi_1^{n_1-n_2}z^{n_1}}{n_1(n_2-t_2)^s}.
\end{eqnarray*}
In other words, this encoding appears to be widespread only as couples of the type\goodbreak$\bxi=(1,1,\ldots,1,\xi_r)$ :
with $\xi_1=1$ and $\omega_1=(1-z)^{-1}{dz}$,
\begin{eqnarray*}
&&\alpha_0^1\left(x_0^{s_1-1}(t_1x_0)^{*s_1}x_1\ldots x_0^{s_{r-1}-1}(t_{r-1}x_0)^{*s_{r-1}}x_{r-1}x_0^{s_r-1}(t_rx_0)^{*s_r}x_r\right)\cr
&&= 
\sum_{n_1>\ldots>n_r}\dfrac{\xi_r^{n_r}}{(n_1-t_1)^{s_1}\ldots(n_r-t_r)^{s_r}}.
\end{eqnarray*}
\end{rema}

\subsection{Duffle relations}\label{ss:stuffl}


Let $\lambda=(\lambda_n)$ be a set of parameters, $\bs=(s_1,\ldots,s_r)$ a composition, $\bxi\in\C^r$. Then
\begin{eqnarray}
\forall n\in\Z_{>0},\quad
M_{\bs,\bxi}^n(\lambda)=\sum_{n>n_1>\ldots>n_r>0}\prod_{i=1}^r\xi_i^{n_i}\lambda_{n_i}^{s_i}
&\mbox{and}&M^n_{(),()}(\lambda)=1.
\end{eqnarray}
We can export the duffle over the tuples $\bs=(s_1,\ldots,s_r)\in\Z_{>0}^r$ and  $\bxi\in\C^r$ with :
\begin{eqnarray}
&&(\bs,\bxi)\duffle((),1)=((),1)\duffle(\bs,\bxi)=(\bs,\bxi)\qquad\quad\mbox{and}\cr
&&(s_1,\bs;\xi_1,\bxi)\duffle(r_1,\br;\rho_1,\brho)\cr
&&=(s_1;\xi_1).\left((\bs;,\bxi)\duffle(r_1,\br;\rho_1,\brho)\right)
+(r_1;\rho_1).\left((s_1,\bs;\rho_1,\bxi)\duffle(\br;\brho)\right)\cr
&&\quad+(s_1+r_1;\xi_1\rho_1).\left((\bs;\bxi)\duffle(\br;\brho)\right)
\end{eqnarray}
\begin{prop}\label{pro:M}
Let  $\bs=(s_1,\ldots,s_l)$ and $\br=(r_1,\ldots,r_k)$ be two compositions,   $\bxi\in\C^l$,  $\brho\in\C^k$.
Then
\begin{eqnarray*}
\forall n\in\N,\quad
M_{\bs,\bxi}^n(\lambda)\,M_{\br,\brho}^n(\lambda)
=
M_{(\bs,\bxi)\duffle(\br,\brho)}^n(\lambda).
\end{eqnarray*}
\end{prop}
\beginproof 
Put the compositions $\bs'=(s_2,\ldots,s_l)$,  $\br'=(r_2,\ldots,r_k)$, the tuples of complex numbers $\bxi'=(\xi_2,\ldots,\xi_l)$ and  $\brho'=(\rho_2,\ldots,\rho_k)$, then\vskip-6mm
\begin{eqnarray*}
&&M_{\bs,\bxi}^n(\lambda)\,M_{\br,\brho}^n(\lambda)\cr
&&=\sum_{n>n_1,n>{n'_1}}
\xi_1^{n_1}\lambda_{n_1}^{s_1}
\,M_{\bs',\bxi'}^{n_1}(\lambda)\,
{\rho}_1^{{n'}_1}\lambda_{{n'}_1}^{r_1}
\,M_{\br',\brho'}^{{n'}_1}(\lambda)\cr
&&=\sum_{n>n_1}
\xi_1^{n_1}\lambda_{n_1}^{s_1}
\,M_{\bs',\bxi'}^{n_1}(\lambda)
\,M_{\br,\brho}^{n_1}(\lambda)
+\sum_{n>{n'}_1}
{\rho}_1^{{n'}_1}\lambda_{{n'}_1}^{r_1}
\,M_{\bs,\bxi}^{{n'}_1}(\lambda)
\,M_{\br',\brho'}^{{n'}_1}(\lambda)\cr
&&+\sum_{n>m}
(\xi_1\rho_1)^{m}
\lambda_{m}^{s_1+r_1}
\,M_{\bs',\bxi'}^{m}(\lambda)
\,M_{\br',\brho'}^{m}(\lambda).
\end{eqnarray*}
A recurrence ended the demonstration.
\endproof
\begin{theo}\label{Th:stuffle}
Let  $\bs=(s_1,\ldots,s_l)$ and $\br=(r_1,\ldots,r_k)$ be two compositions,  $\bxi$ a $l$-tuple and $\brho$ a $k$-tuple of $\calE$, $\bt=(t,\ldots,t)$ a $l$-tuple and $\bt'=(t,\ldots,t)$ a $k$-tuple, both formed by the same parameter $t$ diagonally.  
Then
\begin{eqnarray*}
 \Di(\bF_{\bxi,\bt};\bs)\Di(\bF_{\bxi',\bt'};\bs')
=\Di(\bF_{\bxi'',(t,\ldots,t)};\bs''),
\end{eqnarray*} 
with $(\bs'';\bxi'')=(\bs;\bxi)\duffle(\bs';\bxi')$. 
\end{theo}
\beginproof
With  $\lambda_n=\dfrac1{n-t}$ for all $n\in\N$,  
\quad
$\displaystyle
M_{\bs,\bxi}^n(\lambda)
=\hskip-3pt\sum_{n>n_1>\ldots>n_r}\prod_{i=1}^r\dfrac{\xi_i^{n_i}}{(n_i-t)^{s_i}}.
$
So $\displaystyle \lim_{n\rightarrow\infty} M_{\bs,\bxi}^n(\lambda)= \Di(\bF_{\bxi,\bt};\bs)$ and 
taking the limit of Proposition \ref{pro:M} gives the result.
\endproof
\begin{exam}
The use of examples \ref{ex:stuffle} and \ref{ex:smuffle} gives
\begin{eqnarray*}
&&\Di(\bF_{(\frac23,-1),\bt};(3,1))
\Di(\bF_{(\frac12),(t)};(2))\cr
&&=
\Di(\bF_{(\frac23,-1,\frac12),(t,t,t)};(3,1,2))
+\Di(\bF_{(\frac23,\frac12,-1),(t,t,t)};(3,2,1))\cr
&&+\Di(\bF_{(\frac23,-\frac12),\bt};(3,3))
+\Di(\bF_{(\frac12,\frac23,-1),(t,t,t)};(2,3,1))+\Di(\bF_{(\frac13,-1),\bt};(5,1))
\end{eqnarray*}
\end{exam}
\begin{rema}
Extend the duffle product to triplets $(\bs,\bt,\bxi)\in\cup_{r\in\N^*}\N^r\times\{t\}^r\times\C^r$ by 
\begin{eqnarray*}
(s_1,\bs;t,\bt;\xi_1,\bxi)\duffle(r_1,\br;t,\bt';\rho_1,\brho)
&=&(s_1;t;\xi_1).\left((\bs;\bt;\bxi)\duffle(r_1,\br;t,\bt';\rho_1,\brho)\right)\cr
&+&(r_1;t;\rho_1).\left((s_1,\bs;t,\bt;\rho_1,\bxi)\duffle(\br;\bt';\brho)\right)\cr
&+&(s_1+r_1;t;\xi_1\rho_1).\left((\bs;\bt;\bxi)\duffle(\br;\bt';\brho)\right),
\end{eqnarray*}
and define the function $F$ over $\calI=\cup_{r\in\N^*}\N^r\times\{t\}^r\times\C^r$ by $F(\bs,\bt,\xi)=\Di(\bF_{\bxi,\bt};\bs)$.\\
Then, by Theorem \ref{Th:stuffle}, the function $F:(\calI,\duffle)\rightarrow(\C,.)$ is  morphism of algebras.
\end{rema}

\end{document}